\documentclass[11pt]{amsart}
\usepackage{verbatim}

\theoremstyle{plain}
\newtheorem{theorem}{Theorem}[section]
\newtheorem{proposition}{Proposition}[section]
\newtheorem{corollary}{Corollary}[section]

\newtheorem{conjecture}{Conjecture}[section]

\theoremstyle{definition}
\newtheorem{definition}{Definition}[section]
\newtheorem{example}{Example}[section]
\newtheorem{remark}{Remark}[section]

\newcommand{\s}{{\mathbb S}}
\renewcommand{\r}{{\mathbb R}}
\newcommand{\n}{{\mathbb N}}

\newcommand{\stn}{\s^{n}({\scriptstyle{\frac{1}{\sqrt{2}}}})}
\newcommand{\stm}{\s^{m}({\scriptstyle{\frac{1}{\sqrt{2}}}})}

\newcommand{\ptb}{{\scriptstyle{\left\{\frac{1}{\sqrt{2}}\right\}}}}
\newcommand{\pt}{{\scriptstyle{\frac{1}{\sqrt{2}}}}}
\newcommand{\svm}{\s^{m}({\scriptstyle{\sqrt{\frac{m+1}{m}}}})}
\newcommand{\stl}{\s^{l}({\scriptstyle{\frac{1}{2}}})}

\DeclareMathOperator{\tr}{trace}
\DeclareMathOperator{\grad}{grad}
\DeclareMathOperator{\ricci}{Ricci}
\DeclareMathOperator{\Div}{div}

\begin{document}

\title{The index of biharmonic maps in spheres}

\author{E. Loubeau}
\address{D{\'e}partement de Math{\'e}matiques\\
Laboratoire C.N.R.S. F.R.E. 2218\\
Universit{\'e} de Bretagne Occidentale \\
6, avenue Victor Le Gorgeu\\
France}
\email{loubeau@univ-brest.fr}

\author{C. Oniciuc}
\address{Faculty of Mathematics\\
Al.I. Cuza University of Iasi\\
Bd. Copou Nr. 11\\
6600 Iasi\\
Romania}
\email{oniciucc@uaic.ro}
\curraddr{D{\'e}partement de Math{\'e}matiques\\
Laboratoire C.N.R.S. F.R.E. 2218\\
Universit{\'e} de Bretagne Occidentale \\
6, avenue Victor Le Gorgeu\\
France}

\thanks{The first author thanks the Faculty of Mathematics of Iasi for its hospitality. \\ \indent The second author thanks the C.N.R.S. for support.}
\keywords{Biharmonic maps, stability}
\subjclass{58E20, 53C42}

\begin{abstract}
Biharmonic maps are the critical points of the bienergy functional and generalise harmonic maps. We investigate the index of a class of biharmonic maps, derived from minimal Riemannian immersions into spheres. This study is motivated by three families of examples: the totally geodesic inclusion of spheres, the Veronese map and the Clifford torus.
\end{abstract}

\maketitle

\section{Introduction}
In their original 1964 paper~\cite{E-S1}, Eells and Sampson proposed an infinite dimensional Morse theory on the manifold of smooth maps between Riemannian manifolds. Though their main results concern harmonic maps, they also suggested other functionals. The interest encountered by harmonic maps, and to a lesser account by $p$-harmonic maps, has overshadowed the study of other possibilities, e.g. exponential harmonicity~\cite{E-L}. While the examples cited so far, are all first order functionals, one can investigate problems involving higher derivatives. A prime example of these, is the bienergy, not only for the role it plays in elasticity and hydrodynamics, but also because it can be seen as the next stage, should the theory of harmonic maps fail. 
\newline Witness the case of the two-torus ${\mathbb T}^2$ and the two-sphere ${\mathbb S}^2$: Eells and Wood showed in~\cite{E-W}, that there exists no harmonic map from  ${\mathbb T}^2$ to ${\mathbb S}^2$ (whatever the metrics chosen) in the homotopy classes of Brower degrees $\pm 1$. But the situation is drastically different for biharmonicity. By Palais-Smale, in each and every homotopy class from a surface (or a three-manifold) to a compact manifold, there exists a biharmonic map (cf.~\cite{E-S2}). While this biharmonic map could very well turn out to be also harmonic, this cannot happen to maps of degrees $\pm 1$ between ${\mathbb T}^2$ and ${\mathbb S}^2$.
\newline Let $\phi : (M,g) \to (N,h)$ be a smooth map between Riemannian
manifolds. We define its {\em tension field} to be: 
$\tau(\phi) = \tr{\nabla d\phi}$, 
and, for a compact domain $K \subseteq M$, its {\em bienergy}:
$$E^{2}(\phi) = \frac{1}{2} \int_{K} |\tau(\phi)|^2 \, v_{g} .$$
The critical points of $E^2$, with respect to continuous deformations,
are called {\em biharmonic maps}. The Euler-Lagrange operator associated
to $E^2$ is:
$$ \tau^{2}(\phi) = - \Delta^{\phi} (\tau(\phi)) - \tr{R^{N}(d\phi,
\tau(\phi))d\phi} ,$$
where
$$ \Delta^{\phi} = - \tr_{g}{(\nabla^{\phi})^{2}}
= - \tr_{g}{\left(\nabla^{\phi}\nabla^{\phi} - \nabla^{\phi}_{\nabla}\right)}$$
is the Laplacian on sections of the pull-back bundle $\phi^{-1}TN$ and
$R^N$ is the curvature operator:
$$ R^{N}(X,Y)Z = [\nabla^{N}_{X},\nabla^{N}_{Y}]Z - \nabla^{N}_{[X,Y]}Z , \quad \forall X,Y,Z \in C(TN).$$
Hence, a map is biharmonic if and only if $\tau^{2}(\phi) = 0$.

It is interesting to observe that $\tau^{2}(\phi) = -
J^{\phi}(\tau(\phi))$, where the Jacobi operator $J^{\phi}$ plays an important role in the theory of harmonic maps, since it gives the second variation of the
energy functional $E(\phi) = \frac{1}{2} \int_{M} |d\phi|^2 \, v_{g}$ at
its critical points, the well-known {\em harmonic maps}.
Moreover, the Euler-Lagrange operator related to $E$ is precisely the tension
field $\tau(\phi)$.

The Jacobi operator being linear, harmonic maps
are trivially biharmonic, in fact, global minima of $E^2$. Besides, by means of a Bochner technique, Jiang established a partial converse:
\begin{theorem}~\cite{Jia}\label{th1}
A biharmonic map from a compact domain into a non-positively curved manifold, is harmonic.
\end{theorem}

One can easily construct a non-harmonic biharmonic map, by choosing a
third order polynomial mapping between Euclidean spaces, since, in this
situation, the biharmonic operator is nothing but the Laplacian composed with
itself.

However, starting from a compact manifold, one must choose positively
curved targets, e.g. the sphere. Remembering that, on the one hand,
mapping $\s^n$ into $\s^{n+1}$ as the equator yields a harmonic map,
i.e. a biharmonic map of vanishing bienergy, and that, on the other
hand, sending the whole of $\s^n$ onto the north pole of $\s^{n+1}$ also
has zero bienergy, one can infer the existence of a non-harmonic
biharmonic map from $\s^n$ into $\s^{n+1}$. 
\newline A more precise argument (cf.~\cite{CMO1}) shows that it actually is the embedding of $\stn$ as the tropic, cancer or capricorn, of
$\s^{n+1}$.

In~\cite{CMP}, Caddeo, Montaldo and Piu show that a biharmonic curve on a surface of non-positive Gaussian curvature is a geodesic, i.e. is harmonic, and give examples of non-harmonic biharmonic curves on spheres, ellipses, unduloids and nodoids. On the three-sphere, the only non-harmonic biharmonic curves are certain circles and helices, while $\s^{2}({\scriptstyle{\frac{1}{\sqrt{2}}}})$ is the unique non-harmonic, that is non-minimal, biharmonic surface (\cite{CMO1}). More generally, for the sphere $\s^n$, only circles and helices are
non-harmonic biharmonic curves(\cite{CMO}).

Since biharmonicity derives from a variational principle, one is naturally tempted to take the second variation of the functional $E^2$. While this
was done in general by Jiang(\cite{Jia}), in light of the restrictions
given by Theorem~\ref{th1} on the curvature of the target space, we will
restrict ourselves to maps into the unit sphere. In this configuration,
the second variation formula for the bienergy was obtained by the second author:
 
\begin{theorem}\cite{Oni}
Assume that $(M,g)$ is a compact manifold and $\phi : (M,g) \to
\s^n$ is a biharmonic map. Take $\left\{ \phi_{s,t} \right\}_{s,t
\in\r}$ a two parameter variation of $\phi$ and set
$\frac{\partial \phi_{s,0}}{\partial s}\big|_{s=0} = V$ and 
$\frac{\partial \phi_{0,t}}{\partial t}\big|_{t=0} = W$, then 
$$\frac{\partial^2 E^{2}(\phi_{s,t})}{\partial s \partial t}\big|_{(s,t)=(0,0)}= (I(V),W) = \int_M<I(V),W>v_g,$$
where
\begin{align}\label{eq1}
I(V)&=\Delta^{\phi}(\Delta^{\phi} V)+\Delta^{\phi}\left(\tr{<V,d\phi\cdot>d\phi\cdot}-\vert d\phi\vert^2V\right) \\
& +2<d\tau(\phi),d\phi>V+\vert\tau(\phi)\vert^2V \notag\\
&-2\tr{<V,d\tau(\phi)\cdot>d\phi\cdot}
-2\tr{<\tau(\phi),dV\cdot>d\phi\cdot} \notag\\
& -<\tau(\phi),V>\tau(\phi)+\tr{<d\phi\cdot,\Delta^{\phi} V>d\phi\cdot} \notag\\
&+ \tr{<d\phi\cdot,\tr{<V,d\phi\cdot>d\phi\cdot}>d\phi\cdot}\notag\\
&-2\vert d\phi\vert^2 \tr{<d\phi\cdot,V>d\phi\cdot} \notag\\
&+ 2<dV,d\phi>\tau(\phi)-\vert d\phi\vert^2\Delta^{\phi} V + \vert d\phi\vert^4 V.\notag
\end{align}
\end{theorem}

For this aspect of the theory, the two most important quantities to determine are the nullity and the index:
\begin{definition}
Let $\phi : (M,g) \to \s^n$ be a biharmonic map. The dimension of the
vector space $\left\{ V \in C(\phi^{-1}T\s^n) : I(V) = 0 \right\}$
is called the {\em nullity} of $\phi$.
\end{definition}

The first case to look at, for nullity, is the identity map ${\mathbf 1} : \s^n \to \s^n$, for which the operator $I(V)$ becomes: 
$$ I(V) = \Delta(\Delta V) - 2(n-1) \Delta V + (n-1)^2 V ,$$
and the nullity is the dimension of the vector space $\left\{ V \in
C(T\s^n) : \Delta V = (n-1) V \right\}$. The Hodge decomposition, the
eigenvalues of the sphere and the dimension of the space of Killing
vector fields, show:
\begin{theorem}\cite{Oni}
The identity ${\mathbf 1} : \s^n \to \s^n$ is a biharmonic map of nullity equal to 6, if $n=2$, and $\frac{n(n+1)}{2}$, when $n \geq 3$.
\end{theorem}
\begin{definition}
Let $\phi : (M,g) \to \s^n$ be a biharmonic map. The {\em index} of $\phi$ is the dimension of the largest vector space $\left\{ V\in C(\phi^{-1}T\s^n) : (I(V),V) < 0 \right\}$.
\newline A map of index zero is called {\em stable}, otherwise, it is said to be
{\em unstable}.
\end{definition}

\begin{remark}
i) Any harmonic is clearly a stable biharmonic map.
\newline ii) Equation~\eqref{eq1} shows that $I$ is a linear elliptic self-adjoint differential operator, and the Hilbert space of sections $L^{2}(\phi^{-1}T\s^n)$ is the orthogonal sum of the finite dimensional eigenspaces of $I$. Its spectrum consists of real numbers, bounded from below. Therefore the nullity and the index are finite. 
\end{remark}
 
The purpose of this article is the computation of the index of biharmonic maps
constructed in the following manner:
\begin{theorem}\label{th2}
Let $M$ be a compact manifold, $\psi:M\to
\stn\times\ptb$ a nonconstant harmonic map and ${\mathbf i}: {\stn} \times \left\{ \frac{1}{\sqrt{2}}\right\} \to {\s^{n+1}}$ the inclusion map.  
\newline Then $\phi={\bf i}\circ \psi:M\to
\s^{n+1}$ is nonharmonic biharmonic if and only if
$e(\psi)$ is constant.
\end{theorem}

\begin{proof}
Let $p\in M$ be an arbitrary point, $\{X_i\}_{i=1,\dots,m}$ a
geodesic frame at $p$, and $\eta$ the unit section of the normal bundle. Using the chain rule for the tension field:
\begin{align*}
\tau(\phi)&= d{\bf i}(\tau(\psi))+\tr{\nabla d{\bf
i}(d\psi\cdot,d\psi\cdot)}\\
&=-\sum_{i=1}^{m} <d\psi(X_i),d\psi(X_i)>\eta \\
&=-2e(\psi)\eta,
\end{align*}
so $\phi$ is not harmonic. Here the second fundamental form of the inclusion is 
$\nabla d{\mathbf i} (X,Y) = B(X,Y)= -<X,Y> \eta$.

By a straightforward computation, we obtain, at the point $p$:
\begin{eqnarray*}
\frac{1}{2}\nabla^{\phi}_{X_i}\nabla^{\phi}_{X_i}\tau(\phi)&=&-(X_iX_ie(\psi))\eta + e(\psi)<d\psi(X_i),d\psi(X_i)>\eta \\
&&-2d\phi[(X_ie(\psi))X_i]-e(\psi)\nabla d\psi(X_i,X_i),
\end{eqnarray*}
and
$$
-R^{\s^{n+1}}(d\phi(X_i),\tau(\phi))d\phi(X_i)=
<d\psi(X_i),d\psi(X_i)>\tau(\phi).
$$
Thus
$$
\frac{1}{2}\tau^2(\phi)=(\Delta e(\psi))\eta-2d\phi(\grad{e(\psi)}).
$$
Now, the map $\phi$ is biharmonic if and only if
$$
\Delta e(\psi)=0 \ \hbox{and} \ d\phi(\grad{e(\psi)})=0,
$$
and the theorem follows from the compactness of $M$.
\end{proof}

In the general case we choose to study, we consider a minimal (i.e. harmonic) isometric immersion of a compact Riemannian manifold into the unit sphere and adapt its radius, to obtain $\psi : (M^m ,g) \to \stn$, of constant energy density $\frac{m}{2}$, which is, then, composed with the tropical inclusion into $\s^{n+1}$, in order to attain, according to Theorem~\ref{th2}, a non-harmonic biharmonic map $\phi: ( M^m ,g) \to \s^{n+1}$, of which we will analyse the stability.
\newline This framework will remain throughout the article.

Our general modus operandi will be to look at the pull-back bundle 
$\phi^{-1}T\s^{n+1}$ as three cases:
\begin{enumerate}
\item The one-dimensional {\em normal} sub-bundle, spanned by the vector field
$\eta$ tangent to $\s^{n+1}$ but normal to the tropic $\stn\times\ptb$.
\item The image of the tangent space of the
domain, i.e. $d\phi(TM)$, named the {\em tangent} sub-bundle. 
\item The sub-space consisting of vector fields tangent to the
tropic but orthogonal to the image of the map and referred to as the {\em vertical} sub-bundle.
\end{enumerate}

\section{The normal sub-bundle}

As the image of $\phi$ actually lies in $\stn\times\ptb$, we start by examining vector fields of $C(\phi^{-1}T\s^{n+1})$ of the form $V = f\eta$, $f$ being a function on $M$ and $\eta (p) = (\psi(p),-\pt), \, \forall p\in M$, i.e. normal to the tropical hypersphere.

\begin{proposition}\label{prop1}
Take $V=f\eta$, $f\in C^{\infty}(M)$, then:
$$
(I(V),V)=\int_M  \left(\vert\Delta^{\phi} V-mV\vert^2-4m^2\vert V\vert^2 \right)\, v_g.
$$
Furthermore, if $\Delta f=\lambda f$:
\begin{equation}
(I(V),V)=\left( \lambda^2+4\lambda-4m^2\right)\int_Mf^2 \, v_g.
\end{equation}
\end{proposition}

\begin{proof}
When $V=f\eta$, the various elements of Formula~\eqref{eq1} become:
\begin{align*}
& \tau(\phi)=-m\eta, \quad \tr{<V,d\phi\cdot>d\phi\cdot}=0, \quad \vert d\phi\vert^2=m,\\
& <d\tau(\phi),d\phi>=-m^2, \quad \vert \tau(\phi)\vert^2=m^2, \quad \tr{<V,d\tau(\phi)\cdot>d\phi\cdot}=0, \\
& \tr{<\tau(\phi),dV\cdot>d\phi\cdot}=-md\phi(\grad f), \quad <\tau(\phi),V>\tau(\phi)=m^2V, \\
& \tr{<d\phi\cdot,\Delta^{\phi} V>d\phi\cdot} =(\Delta^{\phi} V)^T, \quad <dV,d\phi>\tau(\phi)=-m^2V,
\end{align*}
and:
$$I(V)=\Delta^{\phi}(\Delta^{\phi} V)-2m\Delta^{\phi} V-3m^2V+2md\phi(\grad f)+ 
(\Delta^{\phi} V)^T .$$
As $V$ is normal to $\stn\times\ptb$:
$$<I(V),V>=<\Delta^{\phi}(\Delta^{\phi} V),V>-2m<\Delta^{\phi} V,V>-3m^2\vert V\vert^2,$$
and integrating by parts yields:
$$(I(V),V)=\int_M\left( \vert \Delta^{\phi} V-mV\vert^2-4m^2\vert V\vert^2 \right) \, v_g.$$
To conclude, we observe that:
$$\Delta^{\phi} V=\Delta^{\phi}(f\eta)=(\Delta f)\eta-2d\phi(\grad f)+ m f\eta .$$
\end{proof}

\begin{remark}
A crucial observation, for this section, is that if $V=f\eta$ and $W=g\eta$, where $f$ and $g$ 
are eigenfunctions orthogonal one to the other, then $(I(V),W) =0$. So counting the eigenfunctions, for which $(I(f\eta),f\eta)$ is negative, together with their multiplicities, computes the index on the normal sub-bundle.
\end{remark}

The principal significance of Proposition~\ref{prop1} is to reveal the influence of the small eigenvalues of $\Delta$ on the index of $\phi$:

\begin{corollary}\label{cor1}
Let $V=f\eta$ and $\Delta f= \lambda f$, then  
$(I(V),V)$ is negative if and only if $\lambda\in [0,2(\sqrt{m^2+1}-1)[$.
\end{corollary}

\begin{example}[Generalised Veronese map]\label{ex3.1}\quad
\newline The well-known Veronese map $\s^{2}({\scriptstyle{\sqrt{3}}}) \to \s^4$ is generalised by the immersion of $\s^{m}({\scriptstyle{\sqrt{\frac{2(m+1)}{m}}}})$ into $\s^{m+p}$, 
$p=\frac{(m-1)(m+2)}{2}$ (cf.~\cite{Che}).
\newline We modify the radius, compose with the tropical inclusion and obtain a non-harmonic biharmonic map $\phi: \svm \to \s^{m+p+1}$. 
\newline As $(\svm , g_{can})$ is isometric to 
$(\s^m , {\scriptstyle{\frac{m+1}{m}}} g_{can})$, its eigenvalues are 
\newline $\left\{ \lambda_{k}= \frac{m}{m+1}k(m+k-1) : k\in \n\right\}$, with
 multiplicity $\frac{(2k+m-1)(k+m-2)!}{k!(m-1)!}$ (\cite{E-R}), and, whatever the dimension, $\lambda_{1} < 2(\sqrt{m^2 +1} -1)$ while $\lambda_{k} > 2(\sqrt{m^2 +1} -1), \, \forall k \geq 2$. 
\newline So, for normal vector fields, contributions to the index of $\phi$, come only from the first two eigenvalues, $\lambda_{0}$ and $\lambda_{1}$, making it greater or equal to $1+(m+1) = m+2$.

Furthermore, the generalised Veronese map being quadratic, it also defines a biharmonic map from ${\r}{\mathrm P}^m$, equipped with the metric $\frac{m+1}{m} g_{can}$, into $\s^{m+p+1}$. But the spectrum, in this case, is reduced to 
$\left\{ \tilde{\lambda}_{k}= \frac{m}{m+1}2k(m+2k-1) : k\in \n\right\}$ and, as  $\tilde{\lambda}_{1}$ is larger than $2(\sqrt{m^2 +1} -1)$, only $\eta$ adds to the index.
\end{example}

\begin{example}[Generalised Clifford torus]\label{ex3.2}\quad
\newline An extension to higher dimensions of the Clifford torus, is the minimal Riemannian immersion $\stl\times\stl \to \s^{m+1}({\scriptstyle{\frac{1}{\sqrt{2}}}}), \, (m=2l)$, which gives rise to a non-harmonic biharmonic map $\phi: \stl\times\stl \to \s^{m+2}$.
\newline The eigenvalues of a product being the sum of the eigenvalues of each factor, the spectrum of the Laplacian on $\stl\times\stl$ is 
\newline $\left\{ \lambda_{k}= 4(p(l+p-1) + q(l+q-1)) : p,q \in \n , \, p+q =k\right\}$. So $\lambda_{1}=4l=2m$ is greater than $2(\sqrt{m^2 +1} -1)$, and, among normal vector fields, only $\eta$ is counted in the index.
\end{example}

Since the spectrum of the Laplacian invariably contains the eigenvalue $0$:

\begin{theorem}
Let $(M^m ,g)$ be a compact Riemannian manifold.
\newline A biharmonic map $\phi: (M^m ,g) \to \s^{n+1}$, obtained as the composition of a minimal isometric immersion $\psi : (M^m ,g) \to \stn$ and the tropical inclusion $\stn\times\ptb \to  \s^{n+1}$, has a strictly positive index, i.e. is unstable.
\end{theorem}

\begin{remark}
This result is clearly reminiscent of the instability of harmonic maps into the sphere~\cite{Leu}.
\end{remark}

To further refine Corollary~\ref{cor1}, one can attempt to control the value of the first non-zero eigenvalue, $\lambda_{1}$, through geometrical constraints. Before anything else, note the upper bound $2m$ on $\lambda_{1}$, imposed by the hypotheses on $\psi$. Indeed, immersing $\stn \overset{j}{\to} \r^{n+1}$ makes each component of the map $j\circ\psi: M \to \r^{n+1}$ an eigenfunction of the Laplacian on $(M,g)$, of eigenvalue $2m$, hence the bound on $\lambda_{1}$.

Given that many examples are built from spheres, a helpful condition is inspired by Einstein manifolds:

\begin{proposition} Let $\psi :(M,g) \to \stn$ be a minimal Riemannian immersion and $\phi : (M,g) \to \s^{n+1}$ the biharmonic map constructed above.
\newline Assume that $\ricci^M \geq \kappa g, \, (\kappa >0)$ and let $\kappa(m) =  \frac{2(m-1)}{m}\left( \sqrt{m^2 +1} - 1 \right)$.
\newline If $\kappa \geq \kappa(m)$, then the contribution of normal vector fields to the index of $\phi$ is 1.
\end{proposition}

\begin{proof} 
All hinges on a theorem of Lichnerowicz (\cite{Lic}), stating that 
$\ricci^M \geq \kappa g$ implies $\lambda_{1}\geq \frac{m}{m-1}\kappa$. 
\newline By Corollary~\ref{cor1}, we only need to ensure that $\frac{m}{m-1}\kappa \geq 2(\sqrt{m^2 +1}-1)$.
\end{proof}

\begin{example}[Totally geodesic inclusion]\label{ex3.3}\quad
\newline A very simple example of a minimal isometric immersion into the sphere, is the totally geodesic inclusion of a lower dimension sphere in the equator, $\s^{m}({\scriptstyle{\frac{1}{\sqrt{2}}}}) \to \stn, \, (m\leq n)$, producing a non-harmonic biharmonic map 
$\phi: \s^{m}({\scriptstyle{\frac{1}{\sqrt{2}}}}) \to \s^{n+1}$. 
\newline The sphere $\s^{m}({\scriptstyle{\frac{1}{\sqrt{2}}}})$ being isometric to $(\s^m, \frac{1}{2} g_{can})$, we can apply the above proposition for $\kappa = 2(m-1)$, and deduce that solely the vector field $\eta$ will contribute to the index of the biharmonic map $\phi$ constructed from this totally geodesic inclusion.
\end{example}

\begin{remark}
Although the projective space of Example~\ref{ex3.1} and the product of spheres of Example~\ref{ex3.2} are Einstein manifolds, and therefore admit a constant $\kappa$ satisfying $\ricci \geq \kappa g$, in both instances $\kappa < \kappa(m)$.
\end{remark}

\section{The tangent sub-bundle}

Vector fields of the tangent sub-bundle are in the image of $TM$ by $d\phi$ and can be written $V = d\phi(X), \, X\in C(TM)$.

\begin{proposition}
Let $V=d\phi(X)$ be a tangential vector field. Then: 
\begin{equation*}
(I(V),V)=\int_M \left( \vert\Delta^{\phi} V+(1-m)V\vert^2-m^2\vert V\vert^2 \right) \, v_g 
\end{equation*}
\end{proposition}

\begin{proof}
As $\phi$ is an isometric immersion and $V$ is in $d\phi(TM)$, and therefore, normal to $\eta$ and $\tau(\phi)$, the components of Equation~\eqref{eq1} are:
\begin{align*}
& \tau(\phi)=-m\eta, \quad  \tr <V,d\phi\cdot>d\phi\cdot=V, \quad \vert d\phi\vert^2=m, \\
& <d\tau(\phi),d\phi>=-m^2,\quad \vert\tau(\phi)\vert^2=m^2,\quad \tr <V,d\tau(\phi)\cdot>d\phi\cdot=-mV,\\
&\tr <\tau(\phi),dV\cdot>d\phi\cdot=mV,\quad <\tau(\phi),V>\tau(\phi)=0,\quad \tr <d\phi\cdot,\Delta^{\phi} V>d\phi\cdot=(\Delta^{\phi} V)^T, \\
&\tr <d\phi\cdot,\tr <V,d\phi\cdot>d\phi\cdot>d\phi\cdot=V, \quad <dV,d\phi>=\Div X ,
\end{align*}
consequently:
$$
I(V)=\Delta^{\phi}(\Delta^{\phi} V)+(1-2m)\Delta^{\phi} V+(\Delta^{\phi} V)^T+(1-2m)V-2m(\Div X)\eta.
$$
Since $V$ is normal to $\eta$ and tangent to the image of $M$:
$$
<I(V),V>=<\Delta^{\phi}(\Delta^{\phi} V),V>+2(1-m)<\Delta^{\phi} V,V>+(1-2m)\vert V\vert^2,
$$
and after an integration by parts:
\begin{equation}\label{eq3}
(I(V),V)=\int_M \left( \vert \Delta^{\phi} V+(1-m)V\vert^2-m^2\vert V\vert^2 \right) \, v_g .
\end{equation}
\end{proof}

Exploiting Formula~\eqref{eq3} requires a kind of Bochner formula linking the Laplacian of the pull-back bundle to the Hodge-de Rham Laplacian on vector fields tangent to $M$.

\begin{proposition}\label{prop4.2}
Let $X$ be a vector field tangent to $M$ and $V=d\phi(X)$ a vector field of $\phi^{-1}T\s^{n+1}$. Then:
\begin{align*}
(I(V),V)& =\int_M \Big( 4\big|\tr{\nabla d\psi(\nabla_{.}X,.)} \big|^2 + 4(\Div X)^2 \\
&+ \vert\Delta_H(X)-2\ricci (X)+ mX\vert^2-m^2\vert X\vert^2 \Big) \, v_g
\end{align*}
\end{proposition}
\begin{proof}
Choosing a geodesic frame $\{X_i\}_{i=1,\dots,m}$ around $p\in M$, means that, at $p$, $\nabla_{X_{i}}X_{j} = 0, \, \forall i,j = 1,\dots,m$, and:
$$\Delta^{\phi} V =-\sum_{i=1}^{m} \nabla^{\phi}_{X_i}\nabla^{\phi}_{X_i}V.$$
Moreover, the second fundamental form of the tropical injection ${\mathbf i}$ is $B(X,Y) = - <X,Y> \eta$, so:
$$\nabla^{\phi}_{X_i}V = \nabla_{X_i}^{\psi}d\psi(X)-<X_i,X>\eta.$$
Composing the preceding expressions implies:
\begin{align*}
\nabla^{\phi}_{X_i}\nabla^{\phi}_{X_i}V&=
\nabla^{\phi}_{X_i}\nabla_{X_i}^{\psi}d\psi(X)
-<X_i,\nabla_{X_i}X>\eta-<X_i,X>d\phi(X_i) \\
&= \nabla d{\mathbf i}(d\psi(X_i),\nabla_{X_i}^{\psi}d\psi(X))
+\nabla_{X_i}^{\psi}\nabla_{X_i}^{\psi}d\psi(X) \\
&-<X_i,\nabla_{X_i}X>\eta-<X_i,X>d\phi(X_i) \\
&=\nabla_{X_i}^{\psi}\nabla_{X_i}^{\psi}d\psi(X)
-<d\psi(X_i),\nabla d\psi(X_i,X)+d\psi(\nabla_{X_i}X)>\eta \\
&-<X_i,\nabla_{X_i}X>\eta-<X_i,X>d\phi(X_i) \\
&=\nabla_{X_i}^{\psi}\nabla_{X_i}^{\psi}d\psi(X)-2<X_i,\nabla
_{X_i}X>\eta
-<X_i,X>d\phi(X_i).
\end{align*}
We conclude that:
\begin{equation}
\label{eq4}
\Delta^{\phi} V=\Delta^{\psi}V+2(\Div X)\eta+V.
\end{equation}
The second half of the proof consists of putting to use the harmonicity of $\psi$ and the curvature tensor of $M$, to express $\Delta^{\psi}V$ in terms of $\nabla d\psi$:
\begin{align*}
\Delta^{\psi}V&=-\sum_{i=1}^{m}\nabla_{X_i}^{\psi}\nabla_{X_i}^{\psi}V \\
&=-\sum_{i=1}^{m} \left( \nabla_{X_i}^{\psi}[\nabla d\psi(X_i,X)]+\nabla d\psi(X_i,\nabla_{X_i}X)+d\psi(\nabla_{X_i}\nabla_{X_i}X) \right),
\end{align*}
because $\nabla^{\psi}_{X_{i}}V =\nabla d\psi(X_i,X)+d\psi(\nabla_{X_i}X)$.

But:
\begin{align*}
\nabla_{X_i}^{\psi}[\nabla d\psi(X_i,X)] &=
\nabla_{X_i}^{\psi}[\nabla d\psi(X,X_i)] \\
&=\nabla_{X_i}^{\psi}\nabla_X^{\psi}d\psi(X_i)
-\nabla_{X_i}^{\psi}d\psi(\nabla_X X_i) \\
&=\nabla_X^{\psi}\nabla_{X_i}^{\psi}d\psi(X_i)
+\nabla_{[X_i,X]}^{\psi}d\psi(X_i) +R^{\psi}(X_i,X)d\psi(X_i) \\
&-\nabla d\psi(X_i,\nabla_X X_{i})- 
d\psi(\nabla_{X_i}\nabla_X X_{i}) \\
&=\nabla_X^{\psi}\nabla_{X_i}^{\psi}d\psi(X_i)
+\nabla_{[X_i,X]}^{\psi}d\psi(X_i)+R^{\psi}(X_i,X)d\psi(X_i) \\ &-d\psi(\nabla_X\nabla_{X_i}X_i + \nabla_{[X_i,X]}X_i + R(X_i,X)X_i) \\
&=\nabla_X^{\psi}\nabla_{X_i}^{\psi}d\psi(X_i)
+\nabla_{[X_i,X]}^{\psi}d\psi(X_i) +2\{<X_i,X>d\psi(X_i)-d\psi(X)\} \\
&-\nabla_X^{\psi}d\psi(\nabla_{X_i}X_i)+\nabla d\psi(X,\nabla_{X_i}X_i) -d\psi(\nabla_{[X_i,X]}X_i)-d\psi(R(X_i,X)X_i) \\
&=\nabla_X^{\psi}[\nabla d\psi(X_i,X_i)] +\nabla d\psi([X_i,X],X_i)\\ &+2\Big(<X_i,X>d\psi(X_i)-d\psi(X)\Big) -d\psi(R(X_i,X)X_i).
\end{align*}
Since, $[X_i,X]=\nabla_{X_i}X$ at $p$:
\begin{align}
\Delta^{\psi}V& = -2\sum_{i=1}^{m}\nabla d\psi(\nabla _{X_i}X,X_i) + d\psi\Big(2(m-1)X -\tr\nabla^2X-\ricci(X)\Big) 
\notag \\
&=-2\sum_{i=1}^{m}\nabla d\psi(\nabla_{X_i}X,X_i)+d\psi\Big (2(m-1)X+\Delta_H(X) -2\ricci(X)\Big)\label{eq5} .
\end{align}
The proposition, then, follows from Equations~\eqref{eq3},~\eqref{eq4} and \eqref{eq5}.
\end{proof}

\begin{theorem}\label{th3}
Let $(M^{m},g)$ ($m\geq 2$) be a compact Riemannian manifold and $\phi = {\mathbf i}\circ \psi : M \to \s^{n+1}$ a biharmonic map, forged as the composition of a minimal isometric immersion $\psi$ and the tropical inclusion $\mathbf i$. 
\newline Let $V$ be a vector field of the pull-back bundle $\phi^{-1}T\s^{n+1}$, such that $V= d\phi(X), \, X\in C(TM)$. 
\newline If $\dim{M} \leq 4$ or $\Div(X) =0$, then $(I(V),V)\geq 0$.
\end{theorem}

\begin{proof}
Let $V=d\psi(X)$, $X\in C(TM)$, via the Yano formula (\cite{Yan}):
$$\int_{M} (\Div{X})^2 - < \tr{\nabla^{2} X} +\ricci{(X)},X> \, v_{g} = \int_{M} \frac{1}{2} |L_{X}g|^2 \, v_{g} ,$$
we know that:
\begin{equation}\label{eq6}
(I(V),V)\geq\int_M \Big( |\tr \nabla^2 X+\ricci(X)|^2 
+m \vert L_Xg\vert^2+2(2-m)(\Div X)^2 \Big)\, v_g. 
\end{equation}
Recall that for $X$ tangent to $(M,g)$:
\begin{equation*} 
\vert L_Xg\vert^2\geq\frac{4}{m}(\Div X)^2 ,
\end{equation*}
since for a local orthonormal frame field $\{X_i\}_{i=1,\dots,m}$ on $M$:
\begin{align*}
\vert L_Xg\vert^2&=\sum_{i,j=1}^{m}\left( (L_Xg)(X_i,X_j)\right)^2
=\sum_{i,j=1}^{m}\left(<\nabla_{X_i}X,X_j>+<\nabla_{X_j}X,X_i>\right)^2\\
&\geq 4\sum_{i=1}^{m}(<\nabla_{X_i}X,X_i>)^2,
\end{align*}
and by the Cauchy inequality:

$$\vert L_Xg\vert^2\geq\frac{4}{m}\left(\sum_{i=1}^{m}<\nabla_{X_i}X,X_i>\right)^2
=\frac{4}{m}(\Div X)^2 .$$
Substituting in \eqref{eq6}, we reach:
$$(I(V),V)\geq 2(4-m)\int_M(\Div X)^2 \, v_g.$$
\end{proof}

\begin{remark}\label{rk3.1}
Using the fact that for a Killing vector field $X$: 
$$\sum_{i=1}^{m}\nabla d\psi(\nabla_{X_i}X,X_i)=0,$$
and $\tr{\nabla^{2} X} + \ricci{X} =0$, one can also show that $I(V)=0$. Moreover, the converse also holds, if $m\leq 4$ or $\Div{X}=0$, since, in these cases, $I(V)=0$ implies $L_{X}g=0$, by Equation~\eqref{eq6}, i.e. $X$ is a Killing vector field.
\end{remark}

\begin{corollary}
The nullity of the biharmonic map $\phi$ is bounded from below by the dimension of ${\mathrm{Isom}}(M,g)$.
\end{corollary}

We can complement Theorem~\ref{th3} when the first eigenvalue of the Laplacian is large enough:
 
\begin{corollary}
Let $\phi : M^{m} \to \s^{n+1}$ be a biharmonic map, fashioned as before, and $\lambda_1$ the first non-zero eigenvalue of the Laplacian on $(M,g)$. 
\newline If $\lambda_1\geq\frac{m^2}{4}$ then $(I(V),V)\geq 0$ for all vector fields $V \in C(\phi^{-1}T\s^{n+1})$ such that $V = d\phi(\grad{f}), \, \Delta f = \lambda f$.
\end{corollary}

\begin{proof}
Let $V=d\phi(X)$ and $X= \grad{f}$, where $\Delta f = \lambda_{1} f$.
\newline The proposition comes from the inequalities:
\begin{align*}
(I(V),V)&=
\int_M \Big( 4\big|\tr{\nabla d\psi(\nabla_{.}X,.)}\big|^2 + 4(\Div X)^2 \\
&+\vert\Delta_H(X)- 2\ricci(X) + m X \vert^2 - m^2\vert X|^2 \Big) \, v_g \\
&\geq \int_M \Big( 4(\Div X)^2-m^2\vert X\vert^2 \Big) \, v_g \\
&=\int_M \Big( 4(\Delta f)^2-m^2\vert\grad f\vert^2 \Big) \, v_g \\
&=\lambda_1(4\lambda_1-m^2)\int_M f^2 \, v_g.
\end{align*}
\end{proof}

\begin{example}[Generalised Clifford torus]\label{ex4.1} \quad
\newline The first eigenvalue of $\stl\times\stl$ is $4l$, so, for the biharmonic 
map $\phi : \stl\times\stl \to \s^{m+2}$ ($m=2l$), $(I(V),V)$ will be positive for vector fields $V=d\phi(\grad{f})$, $f$ eigenfunction of the Laplacian, when $2m \geq \frac{m^2}{4}$, that is $m\leq 8$.
\end{example}

Owing to the nature of our construction, we can link the stability of the biharmonic map $\phi : M \to \s^{n+1}$ to the one of the identity.

\begin{theorem}\label{th4}
Let $\psi$ be a minimal isometric immersion from a compact Riemannian manifold $(M,g)$ into $\stn$, $\phi : M \to \s^{n+1}$ the associated biharmonic map and $I$ its second variation operator.
\newline If the identity ${\mathbf 1}: (M,g) \to (M,g)$ is stable, as a harmonic map, then $(I(V),V)$ is positive for any vector field $V$ of $\phi^{-1}T\s^{n+1}$ tangent to $\phi(M)$.
\end{theorem}

\begin{proof}
Writing a vector field $V\in C(\phi^{-1}T\s^{n+1})$, as $d\phi (X), \, X \in C(TM)$, yields:
$$(I(V),V) = \int_{M} \left( |\Delta^{\phi} V + (1-m)V|^2 -m^2 |V|^2 \right) \, v_{g} .$$
with $\Delta^{\phi} V = \Delta^{\psi} V + 2 (\Div{X}) \eta + V$.
\newline On the other hand, the Jacobi operator of $\psi$ is:
$$ J^{\psi}(V) = \Delta^{\psi}V + \tr{R^{\stn}(d\psi, V) d\psi} ,$$
with
\begin{align*}
\tr{R^{\stn}(d\psi, V) d\psi}&=
\sum_{i=1}^{m} R^{\stn}(d\psi(X_{i}), V) d\psi(X_{i}) \\
&= 2 \sum_{i=1}^{m} \left( <X_{i},X>d\psi(X_{i}) - <X_{i}, X_{i}> V \right) \\
&= 2(1-m)V ,
\end{align*}
so that:
$$J^{\psi}(V) = \Delta^{\psi}V + 2(1-m)V .$$
Therefore:
$$\Delta^{\phi} V = J^{\psi}(V) + 2 (\Div{X}) \eta + (2m-1)V ,$$
and
$$\Delta^{\phi} V + (1-m)V = J^{\psi}(V) + 2 (\Div{X}) \eta + mV .$$
Integrating shows that:
\begin{align*}
(I(V),V) &= \int_{M} \left(|J^{\psi}(V) + mV|^2 + 4(\Div{X})^2 -m^2
|V|^2 \right) \, v_{g} \\ 
&= \int_{M} \left( |J^{\psi}(V)|^2 + 4(\Div{X})^2 + 2m<J^{\psi}(V),V> \right) \, v_{g} \\
&\geq  2m \int_{M} <J^{\psi}(V),V>  \, v_{g} = 2m \int_{M} <J^{\mathbf 1}(V),V>  \, v_{g}.
\end{align*}
\end{proof}

\begin{remark}
The stability of the identity of $(M,g)$ is, in fact, equivalent to $(J^{\psi}(V),V)$ positive for any section $V$ of $\psi^{-1}T\stn$ and tangent to $\phi(M)$.
\newline In this respect, $\phi$ inherits part of the stability of $\psi$.
\end{remark}

\begin{example}[Generalised Clifford torus---Revisited]\label{ex4.2} \quad
\newline The manifold $\stl\times\stl$ being Einstein of constant $\kappa=4(l-1)$, its identity map will be stable if and only if $2\kappa \leq \lambda_{1}$ (\cite{Smi}), i.e. provided $l$ equals $1$ or $2$. For these values of $l$, the index of the biharmonic map $\phi: \stl\times\stl \to \s^{2l+2}$, as we already knew, has no contribution from vector fields of the tangent sub-bundle.
\newline When $l>2$, we need a finer analysis of the sign of $(I(V),V)$.
Staying with a section $V$ of $\phi^{-1}T\s^{2l+2}$, such that $V=d\phi(X)$, $X=\grad{f}$ with $\Delta f= \lambda f$, we know that:
\begin{align}
(I(V),V)&\geq 
\int_{\stl\times\stl}  \Big( 4(\Div X)^2 + \vert\Delta_H(X)- 2\ricci (X) + m X \vert^2 - m^2\vert X|^2 \Big) \, v_g \label{eq7}\\
&= \int_{\stl\times\stl} \Big(4 \lambda^2 f^2 + ( \lambda - 8(l-1) + m)^2 \lambda f^2 - m^2  \lambda f^2 \Big) \, v_{g} \notag\\
&= \lambda \Big( \lambda^2 + 2\lambda(10-3m) + 8m^2 - 48m + 64 \Big) \int_{\stl\times\stl} f^2 \, v_{g} ,\notag
\end{align}
where $m=2l$.
\newline Therefore, when $P(\lambda) = \lambda^2 + 2\lambda(10-3m) + 8m^2 - 48m + 64$ is positive, so is $(I(V),V)$, and, the first non-zero eigenvalue $\lambda_{1}$ of the Laplace-Beltrami operator on $\stl\times\stl$ being $4l=2m$, $P(\lambda_{1}) \geq 0$ if and only if $m \leq 8$. The polynomial $P$ is increasing on the interval $[2m , +\infty [$ when $m\leq 10$, so $P(\lambda_{k}) \geq P(\lambda_{1}) \geq 0, \, \forall k\geq 1$ as long as $m\leq 8$ (this recoups Example~\ref{ex4.1}).
\newline For dimensions superior to $8$, $P(\lambda_{1})$ is negative but $P(\lambda_{2})$, and therefore $P(\lambda_{k})$ for $k\geq 2$, is positive, since $\lambda_{k} = 4p(l+p-1) + 4q(l+q-1) \, (p+q =k)$. 
\newline Nonetheless for $\lambda_{1}$, we can push this investigation a little further, by explicitly computing the term overlooked in Inequality~\eqref{eq7}, i.e. 
$ 4\big|\tr{\nabla d\psi(\nabla_{.}X,.)} \big|^2$. 
\newline Let $f\in C^{\infty}(\stl)$ such that $\Delta^{\stl} f = 2m f$ and $\tilde{f}= f\circ p$, $p$ being the first projection $p: \stl\times\stl \to \stl$. Then $\Delta^{\stl\times\stl} \tilde{f} = 2m \tilde{f}$, $X = \grad{\tilde{f}}= \grad{f}$ and, $f$ being the first eigenvalue of the Laplacian on a sphere, $\nabla_{X_{a}} X = \nabla^{\stl}_{X_{a}} \grad{f} = - 4f X_{a}$, where $\left\{X_{a} \right\}_{a=1,\dots,l}$ is an orthonormal frame of $\stl$.
\newline Then, if $\left\{X_{a}, X_{b} \right\}_{a,b=1,\dots,l}= \left\{X_{i} \right\}_{i=1,\dots,2l}$ is an orthonormal frame of $\stl\times\stl$:
\begin{align*}
4\Big|\tr{\nabla d\psi(\nabla_{.}X,.)} \Big|^2 &= 
4\Big|\sum_{i=1}^{2l}{\nabla d\psi(\nabla_{X_{i}}X,X_{i})} \Big|^2 \\
&= 4\Big|\sum_{a=1}^{l}{\nabla d\psi(\nabla_{X_{a}}X,X_{a})} + 
\sum_{b=1}^{l}{\nabla d\psi(\nabla_{X_{b}}X,X_{b})}\Big|^2 \\
&=  4\Big|-4f\sum_{a=1}^{l}{\nabla d\psi(X_{a},X_{a})} + 
\sum_{b=1}^{l}{\nabla d\psi(0,X_{b})}\Big|^2 \\
&= 64 f^2 \Big|\sum_{a=1}^{l} \sqrt{2} \, \xi \Big|^2 \\
&= 32 m^2 f^2 ,
\end{align*}
and:
$$(I(V),V) = (2m(-8m + 64) + 32 m^2 ) \int_{\stl\times\stl} f^2 \, v_{g} > 0.$$
A similar computation yields the same result for $\tilde{f} = f\circ q$ and $\tilde{f}= f\circ p + f\circ q$, where $q$ is the projection on the second factor.
\newline In conclusion, if $V= d\phi(\grad{f})$ with $\Delta f = \lambda f$, then $(I(V),V) > 0$, whatever the dimension.
\end{example}

\begin{remark}
When our examples have dimension less than four, calling on Theorem~\ref{th3} is a far easier option.
\end{remark}

As in the previous section, when the domain is an Einstein manifold, conditions can be imposed on the Einstein constant to ensure the positivity of $(I(V),V)$.

\begin{proposition}\label{prop4.3}
Let $(M^{m},g)$ be an Einstein manifold of constant $\kappa\in \r$, i.e. 
$\ricci = \kappa g$ and $\phi: M \to \s^{n+1}$ a biharmonic map constructed as previously.
\newline If $\kappa \geq \frac{(m+2)^2}{8}$ then $(I(V),V) \geq 0$, for $V = d\phi(\grad{f}), \,  \Delta f = \lambda f$. 
\end{proposition}

\begin{proof}
Let $X\in C(TM)$ and $V = d\phi(X)$.
\newline Assume that $X=\grad{f}$, $f\in C^{\infty}(M)$ and $\Delta f = \lambda f$.
\newline By Proposition~\ref{prop4.2}:
$$(I(V),V) \geq \int_{M} \Big( 4 (\Div X)^2 + 
|\Delta_{H}(X) - 2\ricci(X) + mX|^2 -m^2 |X|^2 \Big) \, v_{g} ,$$
but $\Delta_{H} X = \Delta_{H}(\grad{f}) = \grad{\Delta f} = \lambda
X$, $\ricci{X} = \kappa X$ and:
\begin{align*}
\int_{M} |X|^2 \, v_{g} &= \int_{M} |\grad{f}|^2 \, v_{g} =
\int_{M} \lambda f^2 \, v_{g} \\
\int_{M} (\Div{X})^2 \, v_{g} &= \int_{M} (\Delta f)^2 \, v_{g} = \lambda^2
\int_{M} f^2 \, v_{g} .
\end{align*}
Therefore:
\begin{align*}
(I(V),V) &\geq \lambda\left( \lambda^2 +2(m+2-2\kappa)\lambda + 4\kappa^2 -4\kappa m )\right) \int_{M} f^2 \, v_{g} ,
\end{align*}
and, if $\kappa \geq \frac{(m+2)^2}{8}$, $\lambda^2 +2(m+2-2\kappa)\lambda + 4\kappa^2 -4\kappa m \geq 0$, so $(I(V),V)$ is positive.
\end{proof}

\begin{example}[Totally geodesic inclusion] \label{ex4.3}\quad
\newline Seeing that the sphere $\s^{m}({\scriptstyle{\frac{1}{\sqrt{2}}}})$ 
has $\kappa = 2(m-1)$ for Einstein constant, we can apply the above proposition to show that the biharmonic map $\phi: \s^{m}({\scriptstyle{\frac{1}{\sqrt{2}}}}) \to \s^{n+1}$ is stable in the direction of vector fields $d\phi(\grad{f})$ ($\Delta f = \lambda f$), as soon as $2(m-1) \geq \frac{(m+2)^2}{8}$, i.e. $m\in [2,9]$.
\newline If $m=1$, then $\kappa =0$ and $(I(V),V)$ is always positive for $V=d\phi(\grad{f})$.
\newline When $m\geq 10$ and $V=d\phi(\grad{f})$ ($\Delta f =\lambda f$):
\begin{align*}
(I(V),V) &= \int_{\s^{m}({\scriptstyle{\frac{1}{\sqrt{2}}}})} 
\Big( 4 (\Div X)^2 + 
|\Delta_{H}(X) - 2\kappa X|^2 + 2m <\Delta_{H}(X) - 2\kappa X,X> \Big) \, v_{g} \\
&= \lambda \big(\lambda^2 + 6(2-m)\lambda + 8(m-1)(m-2) \big) \int_{\s^{m}({\scriptstyle{\frac{1}{\sqrt{2}}}})} f^2 \, v_{g} ,
\end{align*}
which is always positive as the spectrum of the Laplacian on $\s^{m}({\scriptstyle{\frac{1}{\sqrt{2}}}})$ is $\big\{ \lambda_{k} = 2k(m+k-1) : k\in \n \big\}$.
\end{example}

\begin{example}[Generalised Veronese map] \label{ex4.4} \quad
\newline To the minimal Riemannian immersion of $\s^{m}({\scriptstyle{\sqrt{\frac{2(m+1)}{m}}}})$ into $\s^{m+p}$, 
$p=\frac{(m-1)(m+2)}{2}$, Theorem~\ref{th2} associates a non-harmonic biharmonic map from $\phi: \svm \to \s^{m+p+1}$.
\newline Since $\svm$, equipped, as usual, with the canonical metric, is an Einstein manifold of constant $\kappa = \frac{m(m-1)}{m+1}$, we deduce that if $V=d\phi(X)$ for $X=\grad f$, $f\in C^{\infty}(\svm)$ with $\Delta f = \lambda f$, as we saw in the proof of Proposition~\ref{prop4.3}, $(I(V),V)$ is positive as soon as $P(\lambda)= \lambda^2 + 2\lambda(m+2-2\kappa) + 4\kappa^2 - 4m\kappa$ is positive.
\newline Recall that the eigenvalues of the Laplacian on $\svm$ are $\left\{ \lambda_{k} = \frac{m}{m+1}k(m+k-1) : k\in \n\right\}$.
\newline The roots of $P$ are $x_{1} = \frac{m^2 -5m -2}{m+1} - \sqrt{\frac{m^3 -3 m^2 + 16 m + 4}{m+1}}$ and $x_{2} = \frac{m^2 -5m -2}{m+1} + \sqrt{\frac{m^3 -3 m^2 + 16 m + 4}{m+1}}$, and, $\lambda_{k} \geq x_{2}, \, \forall m \geq 2$ and $\forall k\geq 2$, while $\lambda_{1} \geq x_{2}$ only if $m=2,3$ or $4$ (as predicted by Theorem~\ref{th3}). 
\newline If $m\geq 5$, then $\lambda_{1} \in ]x_{1},x_{2}[$ and we need to study in details the term $|\tr{\nabla d\psi(\nabla_{.}X,.)}|^2$, for $X=\grad{f}$ with $\Delta f = \lambda_{1} f$. This is, again, made possible by the fact that, on a sphere, the first eigenfunction of the Laplacian takes the form $f(x) = <u,x>$ with $u \in \r^{m+1}$, and satisfies $\nabla_{X} \grad{f} = -\frac{m}{m+1} f X$.
So, for an orthonormal frame field $\left\{X_{i}\right\}_{i=1,\dots,m}$:
$$ \sum_{i=1}^{m} \nabla d\psi (\nabla_{X_{i}} X , X_{i}) = 
-\frac{m}{m+1} f \sum_{i=1}^{m} \nabla d\psi (X_{i} , X_{i}) = 0,$$
by the harmonicity of $\psi$. 
\newline In conclusion, for the generalised Veronese map, $(I(V),V)$ is positive, when $V=d\phi(\grad{f})$ and $\Delta f = \lambda f$, except for $\lambda=\lambda_{1}$.
\end{example}

\section{The vertical sub-bundle}

The third case of vector fields are sections of the pull-back bundle $\phi^{-1}T\s^{n+1}$, tangent to the tropical sphere $\stn$ but orthogonal to the image of the map.
\newline For such a vector field $V$, the second variation operator is easily worked out:

\begin{proposition}
Let $V\in C(\phi^{-1}T\s^{n+1})$ be orthogonal to $\eta$ and $d\phi(TM)$, then:
\begin{equation}\label{eq8}
(I(V),V) = \int_{M} |\Delta^{\phi} V|^2 -2m <\Delta^{\phi} V , V> \, v_{g}
\end{equation}
\end{proposition}

\begin{proof}
The basic properties of $V$, i.e. $<V,\eta>=<V,d\phi(X)>=0, \, \forall X\in TM$, imply that the non-vanishing constituents of Equation~\eqref{eq1} are:
\begin{align*}
&I(V) = \Delta^{\phi}\Delta^{\phi} V - \Delta^{\phi}(|d\phi|^2 V) + 2 <d\tau(\phi), d\phi> V + 
|\tau(\phi)|^2 V - 2 \tr{<V,d\tau(\phi).>d\phi.} \\
&- 2 \tr{<\tau(\phi), dV.>d\phi.}
+ \tr{< d\phi.,\Delta^{\phi} V> d\phi.} + 2 <dV,d\phi> \tau(\phi) - 
|d\phi|^2 \Delta^{\phi} V + |d\phi|^4 V.
\end{align*}
Taking the inner-product with $V$, produces:
\begin{align*}
<I(V),V> &= <\Delta^{\phi}\Delta^{\phi} V-m\Delta^{\phi} V - 2m^2 V+ m^2 V - 
m \Delta^{\phi} V + m^2 V, V> \\
&= < \Delta^{\phi}\Delta^{\phi} V - 2m \Delta^{\phi} V, V> ,
\end{align*}
and integrating by parts yields:
\begin{equation*}
(I(V),V) = \int_{M} |\Delta^{\phi} V|^2 -2m <\Delta^{\phi} V , V> \, v_{g}
\end{equation*}
\end{proof}

Though Equation~\eqref{eq8} is generally unworkable, assuming that $m=n-1$ and that the normal sub-bundle of $\psi^{-1}T\s^{n}$ is parallelizable enables the expression of $\Delta^{\phi} V$ in function of the shape operator.

\begin{proposition}
Let $\psi: M^m \to \stn$ be a minimal isometric immersion and $\phi: M^m \to \s^{n+1}$ the coupled biharmonic map.
\newline Assume that $m=n-1$ and that there exists a unit section $\xi$ of the normal bundle of $M$ in $\psi^{-1}T\stn$.
\newline Let $V= f\xi, \, f\in C^{\infty}(M)$ and $A_{\xi}$ the shape operator of $\xi$, then:
\begin{align*}
<\Delta^{\psi}V,V> &=  (\Delta f) f + f^2 |\nabla d\psi|^2 ,\\
|\Delta^{\psi} V|^2 &= (\Delta f + |\nabla d\psi|^2 f)^2 +
|2A_{\xi}(d\psi(\grad{f})) + f\tr{(\nabla A_{\xi})(d\psi(.),d\psi(.))}|^2 .
\end{align*}
\end{proposition}

\begin{proof}
Let $V=f\xi$ with $f\in C^{\infty}(M)$, so that $<V,\eta>=<V, d\phi(X)> = 0, \, \forall X \in C(TM)$.
\newline The first thing to observe is that for $X\in C(TM)$:
\begin{align*}
\nabla^{\phi}_{X} V &= \nabla^{\s^{n+1}}_{d\phi(X)} V = \nabla^{\s^{n+1}}_{d\psi(X)} V = \nabla^{\stn}_{d\psi(X)} V - <d\phi(X),V>\eta \\
&= \nabla^{\stn}_{d\psi(X)} V = \nabla^{\psi}_{X} V .
\end{align*}
Furthermore:
\begin{align*}
\nabla^{\phi}_{Y}\nabla^{\phi}_{X} V &= \nabla^{\phi}_{Y}\nabla^{\psi}_{X} V \\
&= \nabla^{\psi}_{Y}\nabla^{\psi}_{X} V - <d\psi(Y), \nabla^{\psi}_{X} V> \eta \\
&= \nabla^{\psi}_{Y}\nabla^{\psi}_{X} V + <\nabla^{\psi}_{X} d\psi(Y),  V> \eta \\
&= \nabla^{\psi}_{Y}\nabla^{\psi}_{X} V + <\nabla d\psi(X,Y) + d\psi(\nabla_{X}Y),  V> \eta \\
&= \nabla^{\psi}_{Y}\nabla^{\psi}_{X} V + <\nabla d\psi(X,Y),  V> \eta ,
\end{align*}
and $\Delta^{\phi} V = \Delta^{\psi} V$, as $\psi$ is harmonic.
\newline Now 
\begin{align*}
\nabla^{\psi}_{X}(f\xi) &= (Xf) \xi + f \nabla^{\psi}_{X} \xi \\
&= (Xf)\xi + f\Big( \nabla^{\perp}_{d\psi(X)} \xi - A_{\xi}(d\psi(X))\Big) . 
\end{align*}
The second order of differentiation is:
\begin{align*}
&\nabla^{\psi}_{Y}\nabla^{\psi}_{X}(f\xi) = 
\nabla^{\psi}_{Y}\left( (Xf)\xi + f\big( \nabla^{\perp}_{d\psi(X)} \xi - A_{\xi}(d\psi(X))\big) \right) \\
&=  Y(Xf)\xi +
(Xf)\big(\nabla^{\perp}_{d\psi(Y)}\xi - A_{\xi}(d\psi(Y))\big) +
\nabla^{\psi}_{Y}\left( f\big( \nabla^{\perp}_{d\psi(X)}\xi - A_{\xi}(d\psi(X)) \big) \right)\\
&= Y(Xf)\xi +
(Xf)\nabla^{\perp}_{d\psi(Y)}\xi - (Xf) A_{\xi}(d\psi(Y)) +
(Yf)\left(\nabla^{\perp}_{d\psi(X)}\xi - A_{\xi}(d\psi(X))\right)\\
+& f\Big( \nabla^{\perp}_{d\psi(Y)}\nabla^{\perp}_{d\psi(X)}\xi -
A_{\nabla^{\perp}_{d\psi(X)}\xi}(d\psi(Y)) -
\nabla_{d\psi(Y)}(A_{\xi}(d\psi(X))) - B(d\psi(Y),A_{\xi}(d\psi(X)))\Big) 
\end{align*}
and
\begin{align*}
\Delta^{\psi}V &= (\Delta f)\xi - 2
\nabla^{\perp}_{d\psi(\grad{f})}\xi + 2 A_{\xi}(d\psi(\grad{f})) \\
 &+ f\left( \Delta^{\perp}\xi + |\nabla d\psi|^2 \xi +
\sum_{i=1}^{m}\left( A_{\nabla^{\perp}_{d\psi(X_{i})}\xi}d\psi(X_{i}) + \nabla_{d\psi(X_{i})}
A_{\xi}(d\psi(X_{i}))\right)\right) \\
=& (\Delta f) \xi + 2 A_{\xi} (d\psi (\grad{f})) 
+ f |\nabla d\psi|^2 \xi 
+ f \sum_{i=1}^{m} \nabla_{d\psi(X_{i})}(A_{\xi}(d\psi(X_{i}))) ,
\end{align*}
since $\nabla^{\perp} \xi = 0$.
\newline As a result:
$$ <\Delta^{\psi}V,V> = (\Delta f) f + f^2 |\nabla d\psi|^2 $$
and:
$$|\Delta^{\psi}V|^2 = (\Delta f + |\nabla d\psi|^2 f)^2 +
|2A_{\xi}(d\psi(\grad{f})) + 
f \sum_{i=1}^{m} \nabla_{d\psi(X_{i})} A_{\xi}(d\psi(X_{i}))\Big|^2 .$$
\end{proof}

\begin{example}[Generalised Clifford torus]\label{ex5.1} \quad 
\newline The minimal isometric immersion of the generalised Clifford torus 
$\stl\times\stl \to \s^{2l+1}({\scriptstyle{\frac{1}{\sqrt{2}}}})$ has codimension one and a parallelizable normal bundle, the vector field $\xi$ being defined by $\xi(p,p') = \sqrt{2}(p,-p')$ for $(p,p') \in \stl\times\stl$. One can readily verify that $|\xi| =1$ and $<\xi , X> =0, \, \forall X\in T\s^{l}\times T\s^{l}$.
\newline Furthermore, the shape operator at $\xi$ is:
\begin{equation*}
A_{\xi} = \sqrt{2} 
\begin{pmatrix}
-\mathrm{Id}  & 0 \\
0 & \mathrm{Id}
\end{pmatrix} ,
\end{equation*}
so, for a function $f$ on $\stl\times\stl$ and an orthonormal frame $\left\{ X_{i}\right\}_{i=1,\dots,2l}$:
\begin{align*}
A_{\xi}(d\psi(\grad{f})) &=\sum_{i=1}^{2l} A_{\xi}(d\psi(X_{i}(f)X_{i}))\\
&= - \sqrt{2} \sum_{i=1}^{l}d\psi(X_{i}(f)X_{i}) + \sqrt{2}
\sum_{i=l+1}^{2l} d\psi (X_{i}(f)X_{i}) ,
\end{align*}
and
$$|A_{\xi}(d\psi(\grad{f}))|^2 = 2 |\grad{f}|^2 .$$
Moreover:
\begin{align*}
\tr{(\nabla A_{\xi})}(d\psi(.),d\psi(.)) &= 
\sum_{i=1}^{2l}(\nabla A_{\xi})(d\psi(X_{i}),d\psi(X_{i})) \\
&= \sum_{i=1}^{2l} (\nabla_{d\psi(X_{i})}A_{\xi}(d\psi(X_{i})) - A_{\xi}(\nabla_{d\psi(X_{i})} d\psi(X_{i}))\\
&= 0 ,
\end{align*}
since $A_{\xi}$ acts as the identity (up to a constant) on each of the factors of $\stl\times\stl$.
\newline The rest follows easily, for a section $V=f\xi$:
\begin{align*}
(I(V),V) &= \int_{\stl\times\stl} |\Delta V|^2 - 2m <\Delta V, V> \, v_{g} \\
=& \int_{\stl\times\stl} (\Delta f + 4l f)^2 + 8 |\grad{f}|^2 - 4 l (f(\Delta f)
+ 4lf^2)  \, v_{g} \\ 
&= \int_{\stl\times\stl} (\Delta f)^2 + 8 l f\Delta f + 16 l^2 f^2 + 8
|\grad{f}|^2 -4 l f(\Delta f) - 16 l^2 f^2  \, v_{g} \\
&= \int_{\stl\times\stl} (\Delta f)^2 + 4l f (\Delta f) + 8 |\grad{f}|^2  \, v_{g} \\
&= \int_{\stl\times\stl} (\Delta f)^2 + 4(l+2) |\grad{f}|^2  \, v_{g} .
\end{align*}
So, for the generalised Clifford torus, $(I(V),V)$ is always positive on the vertical sub-bundle.
\end{example}

\begin{example}[Totally geodesic inclusion]\label{ex5.2} \quad
\newline Note that, though the normal bundle of $\stm$ in $\stn$ is not of codimension one, it is parallelized by $\{e_{m+2},\ldots,e_{n+1}\}$ and we can compute its vertical index.
\newline Let $V=fe_{m+i}, \, (2\leq i\leq n-m)$, $f\in C^{\infty}(\stm)$, then $<V,d\phi(X)>=0, \, \forall X\in C(TM)$ and $<V,\eta>=0$. A short computation shows that $\Delta V=
(\Delta f)e_{m+i}$. 
\newline Assuming that $\Delta f=\lambda f$, we have:
\begin{align*}
(I(V),V)&=\int_{\stm}  \Big( <\Delta V,\Delta V>-2m<\Delta V,V> \Big) \, v_g \\
&=\lambda(\lambda-2m)\int_{\stm} f^2 \, v_g \\
&\geq 0.
\end{align*}
Thus, vector fields $fe_{m+i}$ ($\Delta f=\lambda f$) do not contribute to the index of $\phi : \stm \to \s^{n+1}$.
\end{example}

\section{The index}

For some examples, the study of the three different sub-bundles of the pull-back leads to an estimation of the index:

\begin{proposition}\label{prop6.1}
The index of the non-harmonic biharmonic map $\stm \to \s^{n+1}$, ($m\leq n$), obtained from the totally geodesic inclusion $\stm \to \stn$, is 1.
\newline In particular, the inclusion map $\stn \to \s^{n+1}$ has index 1.
\end{proposition}
\begin{remark}
It is interesting to compare Proposition~\ref{prop6.1} with the index of the totally geodesic embedding $\stm \to \stn$ as a harmonic map, which is $n+1$ if $m\geq 3$ and $n-2$ for $m=2$.
\end{remark}
\begin{proof}
One can easily verify that $(I(V),W)=0$ when $V$ and $W$ are orthogonal vector fields of any of the three cases, except for $V=d\phi(X)$ and $W=d\phi(Y)$, with $\Div{X}=\Div{Y}=0$, and when $V=f\eta$ and $W=d\phi(\grad{f})$ ($\Delta f=\lambda f$, as usual).
\newline In the first situation, a linear combination of $V$ and $W$ can be written as $d\phi(Z)$ with $\Div{Z}=0$, and therefore $\mathrm{span}\left\{V,W\right\}$ does not influence the index.
\newline In the second case, a simple computation leads to:
$$ (I(V),W) = -4\lambda (\lambda +2-2m) \int_{\stm} f^2 \, v_{g}$$
and, since $(I(V),V)$ and $(I(W),W)$ are already known (and positive), we have:
$$(I(V),W)^2 \leq (I(V),V)(I(W),W) ,$$
for any eigenvalue $\lambda$ and dimension $m$, so that the second variation operator is positive on the span of $V$ and $W$.
\newline The index of $\stm \to \s^{n+1}$ is therefore the sum of the different contributions examined in Examples~\ref{ex3.3},\ref{ex4.3} and \ref{ex5.2}.
\end{proof}

\begin{proposition}
The nullity of the non-harmonic biharmonic map $\stm \to \s^{n+1}$, ($m\leq n$) is $\frac{1}{2}(m+1)(m+2) +(m+2)(n-m)$.
\end{proposition}
\begin{proof}
A direct use of the previous computations shows that the operator $I$ preserves the following sub-spaces:
\begin{align*}
i) \, &S_{1}=\left\{ f\eta : f\in C^{\infty}(\stm)\right\} \oplus \left\{ d\phi(\grad{g}) : g\in C^{\infty}(\stm)\right\} \\
ii) \, &S_{2}=\left\{ d\phi(X) : \Div{X}=0 \right\} \\
iii) \, &S_{3}= \left\{f_{1}e_{m+2}+\dots +f_{n-m}e_{n+1} : f_{1},\dots,f_{n-m} \in C^{\infty}(\stm)\right\}
\end{align*}
For the first space, one can easily check that, if $f$ is an eigenfunction, $I(2f\eta + d\phi(\grad{f}))=0$ if and only if $f$ corresponds to the first non-zero eigenvalue $\lambda_{1}=2m$ (which has multiplicity $m+1$). Moreover, any vector field of $S_{1}$ can be written $V=\sum_{i\in \n}2\alpha_{i} f_{i} \eta + \beta_{i} d\phi(\grad{f_{i}})$, with $f_{i}$ eigenfunction corresponding to $\lambda_{i}$. 
\newline If $I(V)=0$, then, by previous remarks, $I(2\alpha_{i} f_{i} \eta + \beta_{i} d\phi(\grad{f_{i}}))=0 \quad \forall i\in \n$, and, necessarily, 
$V= \alpha_{1}(2f_{1}\eta + d\phi(\grad{f_{1}}))$, where $\Delta f_{1}= \lambda_{1} f_{1}$. The contribution of $S_{1}$ to the nullity is $m+1$.
\newline Remark~\ref{rk3.1} shows that the space of Killing vector fields on $\stm$, which has dimension $\frac{m(m+1)}{2}$, is the kernel of $I$ restricted to $S_{2}$.
\newline Finally, for a vertical vector field $V= f_{1}e_{m+2} + \dots + f_{n-m}e_{n+1}$, since $\Delta V= \Delta(f_{1})e_{m+2} + \dots + \Delta(f_{n-m})e_{n+1}$, we see that $I(V)=0$ implies $\Delta\Delta f_{i} -2m \Delta f_{i} =0 \quad \forall i=1,\dots,n-m$, so that each $f_{i}$ is a linear combination of a constant function and eigenfunctions of the first eigenvalue. Consequently, on the vertical sub-bundle, the nullity is $(m+2)(n-m)$.

\end{proof}

However, in general, one cannot hope to have a parallelizable normal sub-bundle and, having to forgo some of the sections, our results can only be lower bounds:

\begin{proposition}
The biharmonic map derived from the generalised Veronese map $\svm \to \s^{m+p}$ ($p=\frac{(m-1)(m+2)}{2}$) has index at least equal to $2m+3$. 
\end{proposition}
\begin{proof}
From the study of Example~\ref{ex3.1}, we know that the index of the generalised Veronese map is at least $m+2$. As we know that $(I(V),V)$ is also negative when $V=d\phi(\grad{f})$, with $\Delta f= \lambda_{1}f$ (and no other eigenvalue), a better lower bound of the index will be obtained if these vector fields combine to span an even larger space on which $(I(V),V)$ is negative definite.
\newline This requires computing $(I(V),W)$ in three different cases: first when $V=d\phi(\grad{f})$ and $W=d\phi(\grad{g})$, $\Delta f =\lambda_{1} f$ and $\Delta g =\lambda_{1} g$ ($f\neq g$); then $V=f\eta$ , $W=d\phi(\grad g)$, $\Delta f = \lambda f$ ($\lambda = \lambda_{0}$ or $\lambda_{1}$) and $\Delta g = \lambda_{1} g$, and, finally, $V=f\eta$ and $W=d\phi(\grad f)$, $\Delta f =\lambda_{1} f$.
\newline If $V=d\phi(\grad{f})$ and $W=d\phi(\grad{g})$, where $f$ and $g$ are two different eigenfunctions of the Laplacian corresponding to the first eigenvalue $\lambda_{1}$ such that $\int_{\svm}  fg \, v_{g} =0$, then:
$$(I(V),W) =  0 .$$
On the other hand, if $V=f\eta$ and $W=d\phi(\grad{g})$, $\Delta f= \lambda f$, $\lambda=\lambda_{0}$ or $\lambda_{1}$, and $\Delta g = \lambda_{1} g$ ($f\neq g$), then:
$$(I(V),W) = 0, $$
since $\int_{\svm} fg \, v_{g} = 0$.
\newline If $f=g$, then $(I(V),W) = \frac{-4m^3}{(m+1)^2} \int_{\svm} f^2 \, v_{g} $ and $(I(V),W)^2 < (I(V),V)(I(W),W)$, so the bilinear form $(I(V),V)$ is definite negative on the space spanned by $f\eta$ and $d\phi(\grad{f})$.
\newline Together, the vector fields $\left\{ \eta, f\eta, d\phi(\grad{g}), \Delta f= \lambda_{1} f, \, \Delta g= \lambda_{1} g\right\}$ (the multiplicity of $\lambda_{1}$ being $m+1$) span a (2m+3)-dimensional space on which $(I(V),V)$ is negative definite.
\end{proof}

Various attempts to find vector fields, other than $\eta$, for which $(I(V),V)$ is negative has led us to:
\begin{conjecture}
The biharmonic map from $\stl\times\stl$ into $\s^{2l+2}$, constructed from the generalised Clifford torus, has index 1.
\end{conjecture}

\end{document}